\documentclass{article}

\RequirePackage[OT1]{fontenc}
\RequirePackage[
amsthm,amsmath
]{imsart}

\usepackage{graphicx}
\usepackage{latexsym,amsmath}
\usepackage{amsmath,amsthm,amscd}
\usepackage{amsfonts}
\usepackage[psamsfonts]{amssymb}
\usepackage{enumerate}

\usepackage{url}

\usepackage{tocvsec2}

\startlocaldefs

\theoremstyle{plain} 
\newtheorem{theorem}{Theorem}

\newtheorem{proposition}[theorem]{Proposition}
\newtheorem*{propzero}{
Carath\'eodory Principle}

\theoremstyle{definition} 
\newtheorem{definition}
{Definition}
\theoremstyle{definition} 

\newtheorem*{ex*}{Example}
\theoremstyle{remark} 

\theoremstyle{remark} 
\newtheorem{remark}
{Remark}
\newtheorem*{remark*}{Remark}
 
\newcommand{\beqa}{\begin{eqnarray}}
\newcommand{\eeqa}{\end{eqnarray}}
 
\newcommand{\bseq}{\begin{subequations}}
 
\newcommand{\eseq}{\end{subequations}}

\newcommand{\dd}{\partial}

\newcommand{\intr}[2]{\overline{#1,#2}}

\renewcommand{\dd}{{\,\operatorname{d}}}

\newcommand{\supp}{\operatorname{supp}}
\newcommand{\card}{\operatorname{card}}

\newcommand{\la}{\lambda}

\newcommand{\de}{\delta}

\newcommand{\Si}{\Sigma}

\newcommand{\R}{\mathbb{R}}

\newcommand{\cc}{{\mathbf{c}}}

\renewcommand{\le}{\leqslant}
\renewcommand{\ge}{\geqslant}

\newcommand{\fl}[1]{\lfloor#1\rfloor}

\endlocaldefs

\begin{document}

\begin{frontmatter}

\title{Tchebycheff systems and extremal problems for generalized moments: \\ a brief survey}
\runtitle{$T$-systems review}

%

\begin{aug}
\author{\fnms{Iosif} \snm{Pinelis}\thanksref{t2}\ead[label=e1]{ipinelis@mtu.edu}}
  \thankstext{t2}{Supported by NSF grant DMS-0805946}
\runauthor{Iosif Pinelis}


\address{Department of Mathematical Sciences\\
Michigan Technological University\\
Houghton, Michigan 49931, USA\\
E-mail: \printead[ipinelis@mtu.edu]{e1}}
\end{aug}

\begin{abstract}
A brief presentation of basics of the theory of Tchebycheff and Markov systems of functions and its applications to extremal problems for integrals of such functions is given. 
The results, as well as all the necessary definitions, are stated in most common terms. 
This work is motivated by specific applications in probability and statistics. 
A few related questions are also briefly discussed, including the one on the existence of a Tchebycheff system on a given topological space. 
\end{abstract}

  
%



\end{frontmatter}

\settocdepth{chapter}


\settocdepth{subsubsection}

\theoremstyle{plain} 


%
The theory of Tchebycheff systems (or, briefly, $T$-systems) presents powerful tools that can be used 
in various extremal problems in analysis, probability, and statistics; see e.g.\ monographs \cite{karlin-studden,krein-nudelman} and the extensive bibliography there. 
However, whereas the Carath\'eodory Principle (as stated below) and the related duality principle (see e.g.\ \cite[Theorem~4.1]{krein-nudelman}) have been used extensively in probability and statistics (see e.g.\ \cite[Chapters XII and XIII]{karlin-studden} and \cite{pin-utev-exp,pin98,normal,winzor}),  
there appear to have been very few (if any) such applications of the theory of $T$-systems, even though the latter offers significant advantages. 
Even in \cite{karlin-studden}, I have not found applications in probability or statistics of the $T$-systems theory per se. 

This brief survey was motivated by specific applications in \cite{radem-asymp}. 
Here we shall present basics of the theory of $T$-systems and its applications to extremal problems for the corresponding generalized moments. 
The results, as well as all the necessary definitions, will be stated in most common terms and thus, it is hoped, 
easy to use; the notions of canonical and principal representations will be avoided here. 
A few related questions will also be briefly discussed, including the one on the existence of a $T$-system on a given topological space. 

\begin{center}***\end{center}

For a nonnegative integer $n$, let $g_0,\dots,g_n$ be (real-valued) continuous functions on a compact 
topological space $X$. 
Let $M$ denote the set of all (nonnegative) Borel measures on $X$. 
Take any point $\cc=(c_0,\dots,c_n)\in\R^{n+1}$ such that 
\begin{equation}\label{eq:M_c}
	M_\cc:=\Big\{\mu\in M\colon\int_X g_i\dd\mu=c_i\quad\text{for all}\quad i\in\intr0n\Big\}\ne\emptyset;  
\end{equation} 

Consider also the condition that the (generalized) polynomial 
\begin{equation}\label{eq:>0}
\text{$\sum_0^n\la_i g_i$ is strictly positive on $X$, for some $(\la_0,\dots,\la_n)\in\R^{n+1}$. }	
\end{equation}

\begin{propzero}
(See e.g.\ \cite{hoeff-extr,karr}.) 
If the topologcal compact space $X$ is Hausdorff and 
\eqref{eq:>0} holds, then the maximum of $\int_X g_{n+1}\dd\mu$ over all $\mu\in M_\cc$ is attained at some measure $\mu_{\max}\in M_\cc$ with $\card\supp\mu_{\max}\le n+1$.  
\end{propzero}

Here, as usual, $\card$ stands for the cardinality and $\supp\mu$ denotes the support set of the measure $\mu$. 
Note that in \cite{karr} it is additionally assumed that $g_0=1$, which is used to provide for the weak compactness of $M_\cc$. However, the same effect is achieved under the more general condition \eqref{eq:>0}.  

\begin{remark}\label{rem:}
The condition that $X$ be compact can oftentimes be circumvented by using, for instance, an appropriate compactification, say $\overline X$, of $X$ if $X$ is (say) only locally compact (as, for instance, $\R^k$ is). 
At that, one may be able to find a function $h$, which is positive and continuous on $X$ and such that the functions $\frac{g_0}h,\dots,\frac{g_n}h$ can be continuously extended from $X$ to $\overline X$; sometimes one of the $g_i$'s can play the role of $h$; or, more generally, $h$ can be constructed as (or based on) a polynomial $\sum_0^n\la_i g_i$. 
Replacing then $g_0,\dots,g_n$ by the continuous extensions of the functions  $\frac{g_0}h,\dots,\frac{g_n}h$ to $\overline X$ and, correspondingly, replacing the measure $\mu$ by the measure $\nu$ (on $X$ and hence on $\overline X$) defined by the formula $\dd\nu=h\dd\mu$, one will largely reduce the original optimization problem on $X$ to one on the compact space $\overline X$. 
\end{remark}

The essential fact is that the upper bound $n+1$ on the cardinality of the support of an extremal measure $\mu$ given in the Carath\'eodory Principle can be approximately halved in the presence of the Tchebycheff or, especially, Markov property. 

\begin{definition}\label{def:T}
The sequence $(g_0,\dots,g_n)$ of functions is a \emph{$T$-system} if the restrictions of these $n+1$ functions to any subset of $X$ of cardinality $n+1$ are linearly independent. 
If, for each $k\in\intr0n$, the initial subsequence $(g_0,\dots,g_k)$ of the sequence $(g_0,\dots,g_n)$ is a $T$-system, then $(g_0,\dots,g_n)$ is said to be an \emph{$M$-system} (where $M$ refers to Markov). 
\end{definition} 
 
By Haar's theorem, linearly independent functions $g_0,\dots,g_n$ on $X$ form a $T$-system on $X$ if and only the problem of best uniform approximation of any given continuous function on $X$ by a polynomial $\sum_0^n\la_i g_i$ has a unique solution; see e.g.\ \cite{schoen}. 

Any $T$-system satisfies the condition \eqref{eq:>0}; see e.g.\ \cite[Theorem~II.1.4]{krein-nudelman}. 

For any $n\ge1$ and any topological space $X$ of cardinality $\ge n+1$, if there exists a $T$-system $(g_0,\dots,g_n)$ of continuous functions on $X$, then $X$ is necessarily Hausdorff. 
Indeed, take any distinct $x_0$ and $x_1$ in $X$. Let $x_2,\dots,x_n$ be any points in $X$ such that $x_0,\dots,x_n$ are distinct. 
The restrictions of the functions $g_0,\dots,g_n$ to the set $\{x_0,\dots,x_n\}$ are linearly independent and hence $g_i(x_0)\ne g_i(x_1)$ for some $i\in\intr0n$. 
Take now any disjoint open sets $O_0$ and $O_1$ in $\R$ containing $g_i(x_0)$ and $g_i(x_1)$, respectively. Then the pre-images $g_i^{-1}(O_0)$ and $g_i^{-1}(O_1)$ of $O_0$ and $O_1$ under the mapping $g_i$ are disjoint open sets in $X$ containing $x_0$ and $x_1$, respectively. Thus, $X$ is a Hausdorff topological space. 

In connection with Remark~\ref{rem:}, it should be noted that, clearly, if $(g_0,\dots,g_k)$ is a $T$-system or an $M$-system, then the same is true of the sequence $(\frac{g_0}h,\dots,\frac{g_n}h)$, for any positive continuous function $h$. 
Note also that, if $(g_0,\dots,g_n)$ is a $T$-system on a set $X'$ containing $X$ and such that $\card(X'\setminus X)\ge n$, then $(h_0,\dots,h_n):=A(g_0,\dots,g_n)$ is an $M$-system on $X$ for some linear (necessarily nonsingular) transformation $A$ of $\R^{n+1}$; cf.\ \cite[Theorem~II.4.1]{krein-nudelman}. 

A $T$-system $(g_0,\dots,g_n)$ with $n\ge1$ on (the compact topological space) $X$ exists only if $X$ is 
one-dimensional (which will be the case in many applications). More precisely, if for some $n\ge1$ there exists a $T$-system of $n+1$ functions on $X$, then $X$ must be homeomorphic to a subset of a circle; for $X\subseteq\R^k$ this was proved in \cite{mairhuber}, and for general $X$ in \cite{curtis} (with an additional restriction) and in \cite{sieklucki}; a further extension of this result to complex $T$-systems was given in \cite{schoen-yang}, where one can also find yet another proof of the real-valued version.  

In fact, the general case of (real-valued) $T$-systems can be easily reduced to the special case with  $X\subseteq\R^k$. Indeed, for any natural $n$ consider the mapping $x\mapsto r(x):=g(x)/\|g(x)\|$ of $X$ into the unit sphere $S^n$ in $\R^{n+1}$, where $g(x):=(g_0(x),\dots,g_n(x))$ and $\|\cdot\|$ is the Euclidean norm. In view of the $T$-property of $(g_0,\dots,g_n)$ and the compactness of $X$, this mapping is correctly defined \big(since $g(x)$ is nonzero for any $x\in X$\big), one-to-one, and continuous, and hence a homeomorphism of the compact Hausdorff set $X$ onto the image in $S^n$ of $X$ under the mapping $r$. In the case $n=1$, this also proves the mentioned result of \cite{mairhuber,curtis,sieklucki}. 
Another elementary observation in this regard, presented in \cite{buck}, is that a $T$-system $(g_0,\dots,g_n)$ with $n\ge1$ on $X$ may exist only if $X$ does not contain a ``tripod'', that is a set homeomorphic to the set $\{(s,0)\in\R^2\colon|s|<1\}\cup\{(0,t)\in\R^2\colon0<t<1\}$.  

We shall henceforth consider the case when $X=[a,b]$ for some $a$ and $b$ such that $-\infty<a<b<\infty$. 
Let $(g_0,\dots,g_n)$ be a $T$-system on $[a,b]$. 
Let $\det\big(g_i(x_j)\big)_0^n$ denote the determinant of the matrix 
$\big(g_i(x_j)\colon i\in\intr0n,\ j\in\intr0n\,\big)$. 
This determinant is continuous in $(x_0,\dots,x_n)$ in the (convex) simplex (say $\Si$) defined by the inequalities $a\le x_0<\dots<x_n\le b$ and does not vanish anywhere on $\Si$. So, $\det\big(g_i(x_j)\big)_0^n$ is constant in sign on $\Si$. 

\begin{definition}\label{def:T+,M+}
The sequence $(g_0,\dots,g_n)$ is said to be a \emph{$T_+$-system} on $[a,b]$ if $\det\big(g_i(x_j)\big)_0^n>0$ for all $(x_0,\dots,x_n)\in\Si$. 
If $(g_0,\dots,g_k)$ is a $T_+$-system on $[a,b]$ for each $k\in\intr0n$, then 
the sequence $(g_0,\dots,g_n)$ is said to be an \emph{$M_+$-system} on $[a,b]$. 
\end{definition}

Clearly, if $(g_0,\dots,g_n)$ is a $T$-system on $[a,b]$, then either $(g_0,\dots,g_n)$ or $(g_0,\dots,g_{n-1},-g_n)$ is a $T_+$-system on $[a,b]$.  
Similarly, if $(g_0,\dots,g_n)$ is an $M$-system on $[a,b]$ then, for some sequence $(s_0,\dots,s_n)\in\{-1,1\}^{n+1}$ of signs, $(s_0g_0\dots,s_n g_n)$ is an $M_+$-system on $[a,b]$. 

In the case when the functions $g_0,\dots,g_n$ are $n$ times differentiable at a point $x\in(a,b)$, consider also the \emph{Wronskians}
\begin{equation*}
	W_0^k(x):=\det\big(g^{(j)}_i(x)\big)_0^k, 
\end{equation*}
where $k\in\intr0n$ and $g^{(j)}_i$ is the $j$th derivative of $g_i$, with $g^{(0)}_i:=g_i$;  
in particular, $W_0^0(x)=g_0(x)$. 

\begin{proposition}\label{prop:W} 
Suppose that the functions $g_0,\dots,g_n$ are (still continuous on $[a,b]$ and) $n$ times differentiable on $(a,b)$. 
Then, for the sequence $(g_0,\dots,g_n)$ to be an $M_+$-system on $[a,b]$, it is necessary that $W_0^k\ge0$ on $(a,b)$ for all $k\in\intr0n$, and it is sufficient that $u_0>0$ on $[a,b]$ and $W_0^k>0$ on $(a,b)$ for all $k\in\intr1n$. 
\end{proposition}

Thus, verifying the $M_+$-property largely reduces to checking the positivity of several functions of only one variable.  

A special case of Proposition~\ref{prop:W} (with $n=1$ and $g_0=1$) is the following well-known fact: if a function $g_1$ is continuous on $[a,b]$ and has a positive derivative on $(a,b)$, then $g_1$ is (strictly) increasing on $[a,b]$; vice versa, if $g_1$ is increasing on $[a,b]$, then the derivative of $g_1$ (if exists) must be nonnegative on $(a,b)$. 

As in this special case, the proof of Proposition~\ref{prop:W} in general can be based on the mean-value theorem; cf.\ e.g.\ \cite[Theorem~1.1 of Chapter~XI]{karlin-studden}, which states that the requirement for $W_0^k$ to be strictly positive on the \emph{closed} interval $[a,b]$ for all $k\in\intr0n$ is equivalent to a condition somewhat stronger than being an $M_+$-system on $[a,b]$; in connection with this, one may also want to look at \cite[Theorem~IV.5.2]{krein-nudelman}. 
Note that, in the applications to the proofs of \cite[Lemmas~2.2 and 2.3]{radem-asymp}, the relevant Wronskians vanish at the left endpoint of the interval. 

The proof of Proposition~\ref{prop:W} can be obtained by induction on $n$ using the recursive formulas for the determinants $\det\big(g_i(x_j)\big)_0^n$ and $W_0^n$ as displayed right above \cite[(5.5) in Chapter VIII]{karlin-studden} and in \cite[(5.6) in Chapter VIII]{karlin-studden}, where we use $g_i$ in place of $\psi_i$. 

\begin{proposition}\label{prop:extr} 
Suppose that $(g_0,\dots,g_{n+1})$ is an $M_+$-system on $[a,b]$ or, more generally, each of the sequences $(g_0,\dots,g_n)$ and $(g_0,\dots,g_{n+1})$ is a $T_+$-system on $[a,b]$. 
Suppose also that condition \eqref{eq:M_c} holds. 
Let $m:=\fl{\frac{n+1}2}$. 
Then one has the following. 
\emph{
\begin{enumerate}[(I)]
	\item \emph{The maximum (respectively, the minimum) of $\int_a^b g_{n+1}\dd\mu$ over all $\mu\in M_\cc$ is attained at a unique measure $\mu_{\max}$ (respectively, $\mu_{\min}$) in $M_\cc$. 
	Moreover, the measures $\mu_{\max}$ and $\mu_{\min}$ do not depend on the choice of $g_{n+1}$, as long as $g_{n+1}$ is such that $(g_0,\dots,g_{n+1})$ is a $T_+$-system on $[a,b]$. 
	}
	\item \emph{There exist subsets $X_{\max}$ and $X_{\min}$ of $[a,b]$ such that $X_{\max}\supseteq\supp\mu_{\max}$, $X_{\min}\supseteq\supp\mu_{\min}$, and 
	}
\begin{enumerate}
	\item \label{even} \emph{if $n=2m$ then $\card X_{\max}=\card X_{\min}=m+1$, $X_{\max}\ni b$, and $X_{\min}\ni a$; }
	\item \label{odd} \emph{if $n=2m-1$ then $\card X_{\max}=m+1$, $\card X_{\min}=m$, and $X_{\max}\supseteq\{a,b\}$. }
\end{enumerate}
\end{enumerate}
}
\end{proposition}

Whereas the maximizer $\mu_{\max}$ and the minimizer $\mu_{\min}$ are each unique in $M_\cc$ for each given $\cc$ with $M_\cc\ne\emptyset$, in particular applications such as the ones in \cite{radem-asymp}, 
one may want to allow the vector $\cc$ to vary, and then the measures $\mu_{\max}$ and $\mu_{\min}$ will vary with $\cc$, and thus the corresponding subsets $X_{\max}$ and $X_{\min}$ of $[a,b]$ may vary. 
Then the number of real variables needed to describe each of the sets $X_{\max}$ and $X_{\min}$ will be about $\frac{n+1}2$, that is, half the number of restrictions on the measure $\mu$ and also half the upper bound on $\card\supp\mu_{\max}$ in the Carath\'eodory Principle; here one should also take into account that, as described in part (II) of Proposition~\ref{prop:extr}, 
the sets $X_{\max}$ and $X_{\min}$ may have to contain at least one of the endpoints $a$ and $b$ of the interval $[a,b]$, with the corresponding reduction in the required number of variables. 
On the other hand, the Carath\'eodory Principle holds for more general systems of functions, defined on a set $X$ of a much more general class. 

To illustrate Proposition~\ref{prop:extr}, one may consider the simplest two special cases when the conditions of the proposition hold and its conclusion is obvious: 
\begin{enumerate}[(i)]
	\item $n=0$, $g_0(x)\equiv1$, $g_1$ is increasing on $[a,b]$, and $c_0\ge0$; then $\supp\mu_{\max}\subseteq\{b\}$ and $\supp\mu_{\min}\subseteq\{a\}$; 
	in fact, $\mu_{\max}=c_0\de_b$ and $\mu_{\min}=c_0\de_a$; 
	here and in what follows, $\de_x$ denotes the Dirac probability measure at point $x$.   
	\item $n=1$, $g_0(x)\equiv1$, $g_1(x)\equiv x$, $g_2$ is strictly convex on $[a,b]$, $c_0\ge0$, and $c_1\in[c_0a,c_0b]$; then 
	$\supp\mu_{\max}\subseteq\{a,b\}$ and $\card\supp\mu_{\min}\le1$; 
	 in fact, $\mu_{\max}=\frac{c_0b-c_1}{b-a}\de_a+\frac{c_1-c_0a}{b-a}\de_b$, and $\mu_{\min}=c_0\de_{c_1/c_0}$ if $c_0>0$ and $\mu_{\min}=0$ if $c_0=0$.  
\end{enumerate}

These examples also show that the $T$-property of systems of functions can be considered as generalized monotonicity/convexity; see e.g.\ \cite{shaked-shanti} and bibliography there. 

\begin{proof}[Proof of Proposition~\ref{prop:extr}] 
Consider two cases, depending on whether $\cc$ is strictly or singularly positive; in equivalent geometric terms, this means, respectively, that $\cc$ belongs to the interior or the boundary of the smallest closed convex cone containing the subset $\{(g_0(x),\dots,g_n(x))\colon x\in[a,b]\}$ of $\R^{n+1}$ \cite[Theorem~IV.6.1]{krein-nudelman}.  

In the first case, when $\cc$ is strictly positive, both statements of Proposition~\ref{prop:extr} follow by \cite[Theorem~IV.1.1]{krein-nudelman}; at that, one should let $X_{\max}=\supp\mu_{\max}$ and $X_{\min}=\supp\mu_{\min}$. (The condition that $\cc$ be strictly positive appears to be missing 
in the statement of the latter theorem; 
cf.\ \cite[Theorem~1.1 of Chapter~1.1]{karlin-studden}.) 

In the other case, when $\cc$ is singularly positive, use \cite[Theorem~III.4.1]{krein-nudelman}, which states that in this case the set $M_\cc$ consists of a single measure (say $\mu_*$), and its support set $X_*:=\supp\mu_*$ is of an index $\le n$; that is, 
$\ell_-+2\ell+\ell_+\le n$, where $\ell_-$, $\ell$, and $\ell_+$ stand for the cardinalities of the intersections of 
$X_*$ with the sets $\{a\}$, $(a,b)$, and $\{b\}$. 
It remains to show that this condition on the index of $X_*$ implies that there exist subsets $X_{\max}$ and $X_{\min}$ of $[a,b]$ satisfying the conditions \eqref{even} and \eqref{odd} of 
Proposition~\ref{prop:extr} and 
such that $X_{\max}\cap X_{\min}\supseteq X_*$. 

If $n=2m$ then $\card(X_*\cap(a,b))=\ell\le\fl{\frac{2m-\ell_--\ell_+}2}
\le\fl{\frac{2m-\ell_-}2}=m-\ell_-$; 
so, $\card(X_*\cup\{b\})\le\ell_-+(m-\ell_-)+1=m+1$. 
Adding now to the set $X_*\cup\{b\}$ any $m+1-\card(X_*\cup\{b\})$ points of the complement of $X_*\cup\{b\}$ to $[a,b]$, 
one obtains a subset $X_{\max}$ of $[a,b]$ such that $X_{\max}\supseteq X_*$, $X_{\max}\ni b$, and $\card X_{\max}=m+1$. 
Similarly, there exists a subset $X_{\min}$ of $[a,b]$ such that $X_{\min}\supseteq X_*$, $X_{\min}\ni a$,  and $\card X_{\min}=m+1$.  

If $n=2m-1$ then $\card(X_*\cap(a,b))=\ell\le\fl{\frac{2m-1-\ell_--\ell_+}2}
\le m-1$ and hence $\card(X_*\cup\{a,b\})\le1+(m-1)+1=m+1$. 
So, there exists a subset $X_{\max}$ of $[a,b]$ such that $X_{\max}\supseteq X_*$, $X_{\max}\supseteq\{a,b\}$, and $\card X_{\max}=m+1$. 
One also has $\card X_*=\ell_-+\ell+\ell_+\le\fl{\frac{2m-1+\ell_-+\ell_+}2}
\le\fl{\frac{2m+1}2}=m$. 
So, there exists a subset $X_{\min}$ of $[a,b]$ such that $X_{\min}\supseteq X_*$ and $\card X_{\min}=m$. 
\end{proof}

\bibliographystyle{abbrv}


\bibliography{C:/Users/Iosif/Documents/mtu_home01-30-10/bib_files/citations}

\end{document}